\documentclass[10pt,twoside,leqno]{siamltex}
\usepackage{amsfonts,epsfig}

\usepackage{MnSymbol}
\usepackage[all]{xy}
\usepackage[active]{srcltx}
\newtheorem{remark}[theorem]{Remark}

\DeclareMathOperator{\im}{Im}

\renewcommand{\ge}{\geqslant}

\begin{document}

%\Large

%  Leave these commented lines here
% \input{elaheader-volx-xx.tex}
% \setcounter{page}{1}

% \renewcommand{\thefootnote}{\fnsymbol{footnote}}
% \renewcommand{\thefootnote}{\arabic{footnote}}
% \renewcommand{\theequation}{\thesection.\arabic{equation}}

\bibliographystyle{plain}
\title{Change of the *congruence canonical
form of 2-by-2 matrices under
perturbations}
% Leave blank; editors will write the exact dates above

\author{
Vyacheslav Futorny\thanks{Department of
Mathematics, University of S\~ao Paulo,
Brazil (futorny@ime.usp.br).}
% Remember to put \and between any two authors
                    \and
Lena Klimenko\thanks{National
Transport University,
 Suvorov 1, Kiev, Ukraine (e.n.klimenko@gmail.com).}
                 \and
Vladimir V.\
Sergeichuk\thanks{Institute of
Mathematics, Tereshchenkivska 3, Kiev,
Ukraine (sergeich@imath.kiev.ua). This
work was done during the visit of V.V.
Sergeichuk to the University of S\~ao
Paulo supported by Fapesp grant
2012/18139-2.}}

% Authors and running title to go on top of each page
\pagestyle{myheadings} \markboth{V.\
Futorny, L.\ Klimenko, and V.V.\
Sergeichuk}{Change of the *congruence
canonical form  of 2-by-2 matrices
under perturbations} \maketitle

\begin{abstract}
It is studied how small perturbations
of a $2\times 2$ complex matrix change
its canonical form for *congruence. The
Hasse diagrams for the closure ordering
on the set of *congruence classes of
$2\times 2$ matrices is constructed.
\end{abstract}

\section{Introduction}
\label{introd}

We study how small perturbations of a
$2\times 2$ complex matrix can change
its *canonical form for *congruence
(matrices $A$ and $B$ are
\emph{*congruent} if $S^*AS=B$ for a
nonsingular $S$). We construct the
closure graph $G_2$, which is defined
for any natural $n$ as follows.

\begin{definition}
The closure graph  $G_n$ for
*congruence classes of $n\times n$
complex matrices is the directed graph,
in which
 each vertex $v$ represents
      in a one-to-one manner a
      *congruence class $C_v$ of
      $n\times n$ matrices, and
  there is a directed path
      from a vertex $v$ to a vertex
      $w$ if and only if one $($and
      hence each$)$ matrix
      from $C_v$ can be transformed
      to a matrix form
      $C_w$ by an arbitrarily small
      perturbation.
\end{definition}

The graph $G_n$ is the Hasse diagram of
the partially ordered set whose
elements are the *congruence classes of
$n\times n$ matrices and $a\preccurlyeq
b$ means that $a$ is contained in the
closure of $b$. Thus, the graph $G_n$
shows how the *congruence classes
relate to each other in the affine
space of $n\times n$ matrices.

The graph $G_n$ is also the closure
graph for *congruence classes of pairs
$(P,Q)$ of $n\times n$ Hermitian
matrices since each $n\times n$ matrix
is uniquely represented in the form
$P+iQ$ and two matrices are *congruent
if and only if the corresponding pairs
are *congruent. Note that the closure
graphs for congruence classes of 2-by-2
and 3-by-3 matrices were constructed in
\cite{f-k-s_cov_congr}.

Unlike perturbations of matrices under
congruence and *congruence,
perturbations of matrices under
similarity and of matrix pencils have
been much studied. The set of Jordan
canonical forms of all matrices that
are arbitrarily close to a given matrix
was described in
\cite{den-thi,mar-par}. This
description was extended to Kronecker's
canonical forms of pencils in
\cite{pok}. The closure graph for
$2\times 3$ matrix pencils was
constructed in \cite{El-ka}. The
article \cite{kag2} develops a
comprehensive theory of closure
relations for similarity classes of
matrices and for equivalence classes of
matrix pencils. The software
StratiGraph \cite{e-j-k} constructs
their closure graphs.

All matrices that we consider are
complex matrices.

\section{The closure graph for
*congruence classes of 2-by-2
matrices}\label{s2}

Define the $n$-by-$n$ matrices:
\begin{equation*}\label{1aa}
J_n(\lambda):=\begin{bmatrix}
\lambda&1&&0\\&\lambda&\ddots&\\&&\ddots&1
\\ 0&&&\lambda
\end{bmatrix},\qquad
\Delta_{n}:=%
\begin{bmatrix}
0 &  &  & 1\\
&  &
\udots & i\\
& 1 &
\udots & \\
1 & i &  & 0
\end{bmatrix}.
\end{equation*}

Each *congruence class contains exactly
one canonical matrix for *congruence;
that allows us to identify the vertices
of $G_n$ with the $n\times n$ canonical
matrices for *congruence. We use the
*congruence canonical matrices given in
\cite{hor-ser_transp} (see also
\cite{hor-ser_can,hor-ser_anyf}).

\begin{lemma}[{\cite{hor-ser_can}}]
\label{le1} Each square complex matrix
is *congruent to a direct sum, uniquely
determined up to permutation of
summands, of matrices of the form
\begin{equation}\label{can}
\begin{bmatrix}0&I_m\\
J_m(\lambda) &0
\end{bmatrix}\
(0\ne \lambda \in\mathbb C,\ |\lambda |<1),
\quad \mu \Delta_{n}\
(\mu \in\mathbb C,\ |\mu |=1) ,\quad
J_k(0).
\end{equation}
\end{lemma}

For each $A\in{\mathbb C}^{n\times n}$,
the set
\begin{equation}\label{msiy}
V_A:=\{C^*A+AC\,|\,C\in{\mathbb
C}^{n\times n}\}
\end{equation}
is a vector space over $\mathbb R$,
which is the tangent space to the
*congruence class of $A$ at the point
$A$ since
\begin{equation}\label{fwt}
(I+\varepsilon
C)^*A(I+\varepsilon C)
=A+\varepsilon(C^*A+AC)
+\varepsilon^2C^*AC
\end{equation}
for all $C\in{\mathbb C}^{n\times n}$
and $\varepsilon\in\mathbb R$. The
numbers $\dim_{\mathbb R} V_A$ and
$n^2-\dim_{\mathbb R} V_A$ are called
the \emph{dimension} and
\emph{codimension over ${\mathbb R}$}
of the *congruence class of $A$.

The following theorem, which is proved
in Section \ref{sss}, is the main
result of the paper.

\begin{theorem}\label{the1}
The closure graph for *congruence
classes of\/ $2\times 2$ matrices is
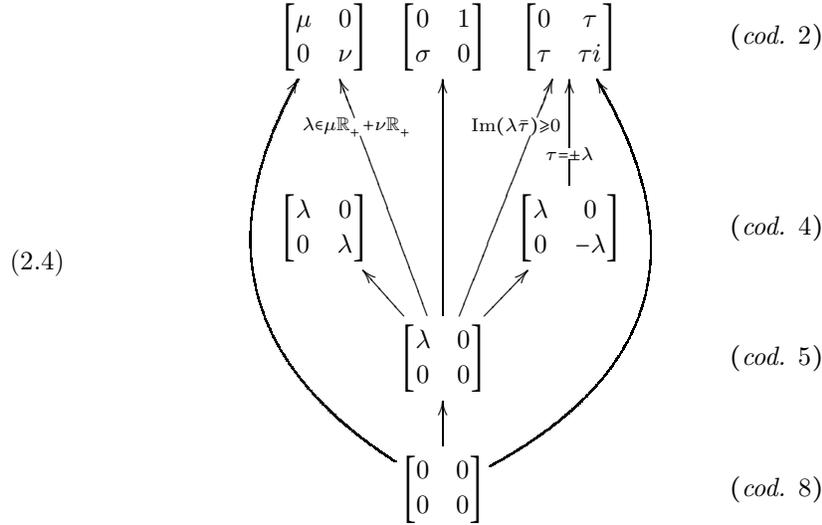
\begin{equation}\label{g1}
\begin{split}\qquad\qquad
\xymatrix@R=17pt@C=8pt{
{\begin{bmatrix} \mu &0\\ 0&\nu
 \end{bmatrix}}
           &
{\begin{bmatrix} 0&1\\ \sigma  &0
\end{bmatrix}}
           &
{\begin{bmatrix} 0&\tau \\
\tau &\tau i \end{bmatrix}}
 &&& {
  (\text{cod. } 2)
}
            %%%%%%%%%%%%%%%%%%%%%%
     \\ \\
            %%%%%%%%%%%%%%%%%%%%%                                   &
{\begin{bmatrix} \lambda &0\\ 0&\lambda
 \end{bmatrix}}
&&
{\begin{bmatrix} \lambda &0\\ 0&-\lambda
 \end{bmatrix}}
\ar[uu]|-(.3){\tau =\pm\lambda}
   &&&
{(\text{cod. }4)}
            %%%%%%%%%%%%%%%%%%%%%%
     \\
            %%%%%%%%%%%%%%%%%%%%%
&{\begin{bmatrix} \lambda &0\\
0&0
 \end{bmatrix}}\ar[ul]\ar[ur]
\ar[uuu]\ar[uuul]|-(.8){\lambda \in
\mu\mathbb R_{_+} +\nu\mathbb
R_{_+}}
\ar[uuur]|-(.8){{\rm Im}
(\lambda\bar\tau)\ge 0\ \quad}& & &&
{(\text{cod. }5)}
            %%%%%%%%%%%%%%%%%%%%%%
     \\
            %%%%%%%%%%%%%%%%%%%%%
&{\begin{bmatrix} 0&0\\
0&0
 \end{bmatrix}}
\ar[u]\ar@/^4.1pc/[uuuul]
\ar@/_4.5pc/[uuuur]
& & && {(\text{cod. }8)}
} \end{split}
\end{equation}
in which $\lambda ,\mu ,\nu ,\sigma
      ,\tau \in\mathbb C$,
      $|\lambda|=|\mu|=
      |\nu|=|\tau|=1$, $|\sigma
      |<1$, $\mathbb R_+$ is the set of
      nonnegative real numbers, and
      $\im (c)$ is the imaginary part
      of $c\in\mathbb C$.
Each *congruence class is given by its
canonical matrix, which is a direct sum
of blocks of the form \eqref{can}. The
graph is infinite: each vertex except
for $\left[\begin{smallmatrix} 0&0\\0&0
\end{smallmatrix}\right]$ represents
an infinite set of vertices indexed by
the parameters of the corresponding
canonical matrix. The *congruence
classes of canonical matrices that are
located at the same horizontal level in
\eqref{g1} have the same codimension
over $\mathbb R$, which is indicated in
round brackets to the right.

The arrow
$\left[\begin{smallmatrix} \lambda &0\\
0&0
\end{smallmatrix}\right]\to
\left[\begin{smallmatrix} \mu &0\\
0&\nu
\end{smallmatrix}\right]$ exists
if and only if $ \lambda
 =\mu a +\nu b$ for some nonnegative
$a,b\in\mathbb R$. The arrow
$\left[\begin{smallmatrix} \lambda &0\\
0&0
\end{smallmatrix}\right]\to
\left[\begin{smallmatrix} 0&\tau \\
\tau &i\tau
\end{smallmatrix}\right]$ exists
if and only if the imaginary part of
$\lambda\bar\tau$ is nonnegative. The
arrow
$\left[\begin{smallmatrix} \lambda &0\\
0&-\lambda
\end{smallmatrix}\right]\to
\left[\begin{smallmatrix} 0&\tau \\
\tau &i\tau
\end{smallmatrix}\right]$ exists
if and only if $\tau =\pm\lambda$. The
arrows
$\left[\begin{smallmatrix} \lambda &0\\
0&0
\end{smallmatrix}\right]\to
\left[\begin{smallmatrix} \lambda &0 \\
0&\pm\lambda
\end{smallmatrix}\right]$ exist
if and only if the value of $\lambda$
is the same in both matrices.  The
other arrows exist for all values of
parameters of their matrices.
\end{theorem}

\begin{remark}\rm
Let $M$ be a $2\times 2$ canonical
matrix for *congruence.

\begin{itemize}
\item If each neighborhood of $M$
    contains a matrix whose
    *congruence canonical form is
$N$, then there is a directed path
    from $M$ to $N$ in \eqref{g1}
    (if $M=N$, then there is always
    the ``lazy'' path of length $0$
    from $M$ to $N$).

  \item The closure of the
      *congruence class of $M$ is
      equal to the union of the
      *congruence classes of all
      canonical matrices $N$ such
      that there is a directed path
      from $N$ to $M$ (if $M=N$
      then there is always the
      ``lazy'' path).
\end{itemize}
\end{remark}

\section{Proof of Theorem \ref{the1}}
\label{sss}

The following lemma is given in
\cite[Example 2.1]{def-sesq}; it is a
special case of \cite[Theorem
2.2]{def-sesq} in which an analogous
statement is given for all square
matrices.

\begin{lemma}
\label{le_def_a} Let $A$ be any
$2\times 2$ matrix. Then all matrices
$A+E$ that are sufficiently close to
$A$ can be simultaneously reduced by
transformations
\begin{equation}\label{dtg}
A+E\mapsto {\cal
S}(E)^* (A+E) {\cal
S}(E),\quad\text{${\cal S}(E)$
is nonsingular and holomorphic
at $0$,}
\end{equation}
to one of the following forms:
\begin{align*}
&\begin{bmatrix} 0&0\\
0&0
 \end{bmatrix}+
\begin{bmatrix}
 *&*\\ *&*
 \end{bmatrix},
                            &&
\begin{bmatrix} \lambda &0\\
0&0
 \end{bmatrix}+
\begin{bmatrix}
\varepsilon_{\lambda}&0\\ *&*
 \end{bmatrix}\ (|\lambda |=1) ,
                         \\&
\begin{bmatrix} \lambda &0\\
0&\pm\lambda
 \end{bmatrix}+
\begin{bmatrix}
\varepsilon_{\lambda}&0\\ *&
\delta_{\lambda}
 \end{bmatrix}\ (|\lambda |=1),
                         &&
\begin{bmatrix} \lambda &0\\
0&\mu
 \end{bmatrix}+
\begin{bmatrix}
 \varepsilon_{\lambda}&0\\
 0&\delta_{\mu}
 \end{bmatrix}
\begin{array}{l}
(\lambda \ne\pm\mu, \\
\ \ |\lambda|=|\mu |= 1),  \\
\end{array}
                       \\&
\begin{bmatrix} 0&1\\
\lambda  &0
 \end{bmatrix}+
\begin{bmatrix}
 0&0\\ *&0
 \end{bmatrix}\ (|\lambda | < 1),
                         &&
\begin{bmatrix} 0&\lambda \\
\lambda  &\lambda i
 \end{bmatrix}+
\begin{bmatrix}
 *&0\\ 0&0
 \end{bmatrix}\ (|\lambda|=1).
 \hspace{-12pt}
\end{align*}
Each of these matrices has the form
$A_{\rm can}+{\cal D}$, in which
$A_{\rm can}$ is a direct sum of blocks
of the form \eqref{can}, the stars in
${\cal D}$ are complex numbers, all
$\varepsilon_{\lambda},\delta_{\lambda},
\delta_{\mu}$ are either real numbers
if $\lambda ,\mu \notin\mathbb R$ or
pure imaginary numbers if $\lambda ,\mu
\in\mathbb R$. $($Clearly, ${\cal D}$
tends to zero as $E$ tends to zero.$)$
Twice the number of stars plus the
number of entries
$\varepsilon_{\lambda},\delta_{\lambda},
\delta_{\mu}$ is the smallest that can
be attained by using transformations
\eqref{dtg}; it is equal to the
codimension over $\mathbb R$ of the
*congruence class of $A$.
\end{lemma}

Note that the codimensions of
congruence and *congruence classes were
calculated in \cite{ter-dor,f-ser} and
\cite{t_d_*,def-sesq}, respectively.

By \cite[Part III, Theorem 1.7]{tan},
the boundary of each *congruence class
is a union of *congruence classes of
strictly lover dimension, which ensures
the following lemma.

\begin{lemma}
\label{lej} If $M\to N$ is an arrow in
the closure graph $G_2$, then the
*congruence class $C_M$ of $M$ is
contained in the closure of the
*congruence class $C_N$ of $N$, and so
the codimension of $C_M$ is greater
than the codimension of $C_N$.
\end{lemma}

The proof of Theorem \ref{the1} is
divided into two steps.

\subsection*{Step 1: Let us prove that
each arrow in \eqref{g1} is correct}

To make sure that an arrow $M\to N$ is
correct, we need to prove that the
canonical matrix $M$ can be transformed
by an arbitrarily small perturbation to
a matrix whose *congruent canonical
form is $N$. Consider each of the
arrows of \eqref{g1}. \medskip

$\bullet$ \emph{The arrows
${\left[\begin{smallmatrix}
0&0\\
0&0
 \end{smallmatrix}\right]}
\to{\left[\begin{smallmatrix} \mu  &0\\
0&\nu
 \end{smallmatrix}\right]}
$, ${\left[\begin{smallmatrix}
0&0\\
0&0
 \end{smallmatrix}\right]}
\to{\left[\begin{smallmatrix} \lambda
 &0\\
0&0
 \end{smallmatrix}\right]}
$, and ${\left[\begin{smallmatrix}
0&0\\
0&0
 \end{smallmatrix}\right]}
\to{\left[\begin{smallmatrix} 0 &\tau \\
\tau &\tau i
 \end{smallmatrix}\right]}
$ are correct.}

Indeed,  each $2\times 2$ matrix $A$ is
*congruent to $\varepsilon A$, in which
$\varepsilon $ is any positive real
number, and each neighborhood of
$\left[\begin{smallmatrix}
0&0\\
0&0
 \end{smallmatrix}\right]$ contains
$\varepsilon A$ with a sufficiently
small $\varepsilon $. We set $A:=
\left[\begin{smallmatrix} \mu  &0 \\
0&\nu i
 \end{smallmatrix}\right]$,
 $\left[\begin{smallmatrix} \lambda  &0 \\
0&0
 \end{smallmatrix}\right]$, or $\left[\begin{smallmatrix} 0 &\tau \\
\tau &\tau i
 \end{smallmatrix}\right]$.

\medskip

$\bullet$ \emph{The arrow
${\left[\begin{smallmatrix}
\lambda   &0\\
0&0
 \end{smallmatrix}\right]}
\to{\left[\begin{smallmatrix} \mu  &0\\
0&\nu
 \end{smallmatrix}\right]}
$ $(|\lambda |=|\mu |=|\nu |=1)$ exists
if and only if $ \lambda \in \mu\mathbb
R_{+} +\nu\mathbb R_{+} =\{\mu a +\nu b
|a,b\in\mathbb R, a\ge 0, b\ge 0\}$
$($in particular,
${\left[\begin{smallmatrix}
\lambda   &0\\
0&0
 \end{smallmatrix}\right]}
\to{\left[\begin{smallmatrix} \lambda  &0\\
0&\lambda
 \end{smallmatrix}\right]}$
and ${\left[\begin{smallmatrix}
\lambda   &0\\
0&0
 \end{smallmatrix}\right]}
\to{\left[\begin{smallmatrix} \lambda  &0\\
0&-\lambda
 \end{smallmatrix}\right]}$
exist$)$.}

Let ${\left[\begin{smallmatrix}
\lambda &0\\
0&0
 \end{smallmatrix}\right]}
\to{\left[\begin{smallmatrix} \mu &0\\
0&\nu
 \end{smallmatrix}\right]}
$ exist; i.e., there exists an
arbitrarily small perturbation
\begin{equation}\label{ltr}
\begin{bmatrix}
\lambda &0\\
0&0
 \end{bmatrix}+E:=\begin{bmatrix}
     \lambda  +\varepsilon _{11}
     & \varepsilon_{12}  \\
     \varepsilon_{21} & \varepsilon _{22}  \\
   \end{bmatrix}\quad\text{of}\quad
 \begin{bmatrix}
\lambda  &0\\
0&0
 \end{bmatrix}
\end{equation}
that is
 *congruent to
${\left[\begin{smallmatrix} \mu &0\\
0&\nu
 \end{smallmatrix}\right]}$.
This means that there exists a
nonsingular $S=\left[\begin{smallmatrix}
x&y\\
z&t
 \end{smallmatrix}\right]
$ such that \[S^*\begin{bmatrix} \mu &0\\
0&\nu
 \end{bmatrix}S=
 \begin{bmatrix}
\lambda   &0\\
0&0
  \end{bmatrix}+E,\]
i.e.,
\begin{equation}\label{hyr}
 \begin{alignedat}{4}
\bar xx\mu +\bar zz\nu &=\lambda
+\varepsilon _{11}
   &\qquad
\bar xy\mu
+\bar zt\nu &=\varepsilon _{12}
     \\
\bar yx\mu+\bar tz\nu &=\varepsilon _{21}
&
    \bar yy\mu +\bar tt\nu &=
    \varepsilon _{22}
      \end{alignedat}
\end{equation}
Clearly, $\lambda\in\mu\mathbb R_+
+\nu\mathbb R_+$ if and only if there
are $x$ and $z$ and an arbitrarily
small $\varepsilon _{11}$ such that the
first equality in \eqref{hyr} holds (we
can take $\varepsilon _{11}=0$). Taking
arbitrarily small $y,t$ for which $S$
is nonsingular, we get arbitrarily
small $\varepsilon _{12},\varepsilon
_{21},\varepsilon _{22}$ for which the
other equalities in \eqref{hyr} hold.
\medskip

$\bullet$ \emph{The arrow
${\left[\begin{smallmatrix}
\lambda  &0\\
0&0
 \end{smallmatrix}\right]}
\to{\left[\begin{smallmatrix} 0 &1\\
\sigma &0
 \end{smallmatrix}\right]}
$ $(|\lambda  |=1$, $|\sigma |<1)$
exists for all $\lambda $ and
$\sigma$.}

Let ${\left[\begin{smallmatrix}
\lambda &0\\
0&0
 \end{smallmatrix}\right]}
\to{\left[\begin{smallmatrix} 0 &1\\
\sigma&0
 \end{smallmatrix}\right]}
$; i.e., there exists an arbitrarily
small perturbation \eqref{ltr} that is
 *congruent to
${\left[\begin{smallmatrix} 0 &1\\
\sigma&0
 \end{smallmatrix}\right]}$.
This means that there exists a
nonsingular
$S=\left[\begin{smallmatrix}
x&y\\
z&t
 \end{smallmatrix}\right]
$ such that
\[S^*\begin{bmatrix}
0 &1\\
\sigma&0
 \end{bmatrix}S=
 \begin{bmatrix}
\lambda &0\\
0&0
  \end{bmatrix}+E,\]
i.e.,
\begin{equation}\label{kue}
      \begin{alignedat}{4}
\bar xz +\bar z\sigma  x &=\lambda
+\varepsilon _{11}
   &\qquad
\bar xt
+\bar z\sigma y &=\varepsilon _{12}
     \\
\bar yz +\bar t\sigma x &=\varepsilon _{21}
 &
    \bar yt +\bar t\sigma y &=\varepsilon _{22}
      \end{alignedat}
    \end{equation}

Let $\bar zx=u+iv$, $\sigma=\alpha
+\beta i$, and $\lambda +\varepsilon
_{11}=a+bi$, in which $u,v,\alpha
,\beta ,a,b\in\mathbb R$, then the
first equation in \eqref{kue} takes the
form $ (u-vi)+(u+vi)(\alpha +\beta
i)=a+bi$, which gives the system
\begin{align*}
(1+\alpha)u-\beta v &=a\\
\beta u +(\alpha-1)v &=b
\end{align*}
with respect to the unknowns $u$ and
$v$. Its determinant $\alpha ^2+\beta
^2-1$ is nonzero since $|\sigma |<1$.
Therefore, the first equation in
\eqref{kue} holds for some $x$ and $z$.
Taking arbitrarily small $y,t$ for
which $S$ is nonsingular, we get
arbitrarily small $\varepsilon
_{12},\varepsilon _{21},\varepsilon
_{22}$ for which the other equalities
in \eqref{kue} hold. \medskip

$\bullet$  \emph{The arrow
${\left[\begin{smallmatrix}
\lambda  &0\\
0&0
 \end{smallmatrix}\right]}
\to{\left[\begin{smallmatrix} 0 &\tau \\
\tau &\tau i
 \end{smallmatrix}\right]}
$ $(|\lambda  |=|\tau |=1)$ exists if
and only if $\im(\lambda\bar\tau )\ge 0
$.}

Let ${\left[\begin{smallmatrix}
\lambda &0\\
0&0
 \end{smallmatrix}\right]}
\to{\tau \left[\begin{smallmatrix} 0 &1\\
1&i
 \end{smallmatrix}\right]}
$ exist; i.e., there exists an
arbitrarily small perturbation
\eqref{ltr} that is
 *congruent to
$\tau {\left[\begin{smallmatrix} 0 &1\\
1&i
 \end{smallmatrix}\right]}$.
This means that there exists a
nonsingular
$S=\left[\begin{smallmatrix}
x&y\\
z&t
 \end{smallmatrix}\right]
$ such that \[S^*\tau \begin{bmatrix}
0 &1\\
1&i
 \end{bmatrix}S=
 \begin{bmatrix}
\lambda &0\\
0&0
  \end{bmatrix}+E,\]
i.e.,
\begin{equation}\label{kmy}
      \begin{alignedat}{4}
\bar zx +\bar xz +\bar zz i&
=\bar\tau (\lambda
+\varepsilon _{11})
   &\qquad
\bar zy
+\bar xt+\bar zti &=\bar\tau \varepsilon _{12}
     \\
\bar tx +\bar yz+\bar tz &=\bar\tau
\varepsilon _{21}
 &
    \bar ty +\bar yt+\bar tti &=\bar\tau
    \varepsilon _{22}
      \end{alignedat}
    \end{equation}

Consider the first equation in
\eqref{kmy}. Since $\bar\tau (\lambda
+\varepsilon _{11})\ne 0$, $z\ne 0$
too. Thus, \[\im(\bar\tau (\lambda
+\varepsilon _{11}))=\im(\bar zx +\bar
xz +\bar zz i)=\bar zz>0,\] and so
$\im(\bar\tau \lambda)\ge 0$.

Conversely, if $\im(\bar\tau
\lambda)\ge 0$, then we can take $x ,z
$ and an arbitrarily small $\varepsilon
_{11}$ such that the first equation in
\eqref{kmy} holds. Taking arbitrarily
small $y,t$ for which $S$ is
nonsingular, we get arbitrarily small
$\varepsilon _{12},\varepsilon
_{21},\varepsilon _{22}$ for which the
other equalities in \eqref{kmy} hold.
\medskip

$\bullet$ \emph{The arrow
${\left[\begin{smallmatrix}
\lambda  &0\\
0&-\lambda
 \end{smallmatrix}\right]}
\to{\left[\begin{smallmatrix} 0 &\tau \\
\tau &\tau i
 \end{smallmatrix}\right]}
$ $(|\lambda  |=|\tau |=1)$ exists if
and only if $\lambda=\pm\tau $.}

Let ${\left[\begin{smallmatrix}
\lambda &0\\
0&-\lambda
 \end{smallmatrix}\right]}
\to{\tau \left[\begin{smallmatrix} 0 &1\\
1&i
 \end{smallmatrix}\right]}
$ exist; i.e., there exists an
arbitrarily small perturbation
$\left[\begin{smallmatrix}
\lambda &0\\
0&-\lambda
 \end{smallmatrix}\right]+E$ of $\left[\begin{smallmatrix}
\lambda &0\\
0&-\lambda
 \end{smallmatrix}\right]$ that is
 *congruent to
$\tau {\left[\begin{smallmatrix} 0 &1\\
1&i
 \end{smallmatrix}\right]}$.
This means that there exists a
nonsingular $S$ such that \[S^*\tau
\begin{bmatrix}
0 &1\\
1&i
 \end{bmatrix}S=
 \begin{bmatrix}
\lambda &0\\
0&-\lambda
  \end{bmatrix}+E.\]
Equating the determinants of both
sides, we find that  $ -\tau ^2\det
(S^*S)$  is arbitrarily close to
$-\lambda ^2$. Since \[\det
(S^*S)=\overline{\det S}\det S
\in\mathbb R\quad\text{and}\quad
|\lambda |=|\tau |=1,\] we have
$\lambda=\pm\tau $.

Conversely, let $\lambda=\pm\tau $.
Since \[ \begin{bmatrix}
  1 & 1 \\
  1/2 & -1/2 \\
 \end{bmatrix}
 \begin{bmatrix}
  1 & 0\\
  0& -1 \\
 \end{bmatrix}
\begin{bmatrix}
  1 & 1/2 \\
  1& -1/2 \\
 \end{bmatrix}=
\begin{bmatrix}
  0& 1 \\
  1 & 0 \\
 \end{bmatrix}
\]
and
\[
\begin{bmatrix}  1 & 0 \\
  0& -1 \\
 \end{bmatrix}
\begin{bmatrix}
  0& 1 \\
  1& 0 \\
 \end{bmatrix}
 \begin{bmatrix}
  1 & 0 \\
  0& -1 \\
 \end{bmatrix}=
\begin{bmatrix}
  0 & -1 \\
  -1& 0\\
 \end{bmatrix},
\]
we have that $\left[\begin{smallmatrix}
1&0\\
0&-1
 \end{smallmatrix}\right]$ is
 *congruent to $\pm
 \left[\begin{smallmatrix}
0&1\\
1&0
 \end{smallmatrix}\right]$.
Its arbitrarily small perturbation $\pm
 \left[\begin{smallmatrix}
0&1\\
1&\varepsilon i
 \end{smallmatrix}\right]$
 ($\varepsilon\in \mathbb R$, $\varepsilon>0$) is
 *congruent to $\pm
 \left[\begin{smallmatrix}
0&1\\
1&i
 \end{smallmatrix}\right]$
via $\diag(\sqrt \varepsilon,1/\sqrt
\varepsilon) $. Therefore,
${\left[\begin{smallmatrix}
1 &0\\
0&-1
 \end{smallmatrix}\right]}
\to{\pm\left[\begin{smallmatrix} 0 &1\\
1&i
 \end{smallmatrix}\right]}
$, and so ${\left[\begin{smallmatrix}
\lambda &0\\
0&-\lambda
 \end{smallmatrix}\right]}
\to{\tau \left[\begin{smallmatrix} 0 &1\\
1&i
 \end{smallmatrix}\right]}
$.

\subsection*{Step 2: Let us prove that we
have not missed arrows in \eqref{g1}}

We write $M\nrightarrow N$ if the
closure graph $G_2$ does not have the
arrow $M\to N$; i.e., each matrix
obtained from $M$ by an arbitrarily
small perturbation is not *congruent to
$N$. The evident statement ``if
$L\leftarrow M\nrightarrow N$ then
$L\nrightarrow N$'' and Lemma \ref{lej}
ensure that we only need to prove the
absence of the arrows
\[{\left[\begin{smallmatrix}
\lambda  &0\\
0&\pm\lambda
 \end{smallmatrix}\right]}\to
{\left[\begin{smallmatrix} \mu  &0\\
0&\nu
 \end{smallmatrix}\right]},\quad
 {\left[\begin{smallmatrix}
\lambda  &0\\
0&\pm\lambda
 \end{smallmatrix}\right]}\to
{\left[\begin{smallmatrix} 0 &1\\
\sigma  &0
 \end{smallmatrix}\right]},\quad
 {\left[\begin{smallmatrix}
\lambda  &0\\
0&\lambda
 \end{smallmatrix}\right]}\to
{\left[\begin{smallmatrix} 0 &\tau \\
\tau &\tau i
 \end{smallmatrix}\right]}.\]

$\bullet$  ${\left[\begin{smallmatrix}
\lambda  &0\\
0&\pm\lambda
 \end{smallmatrix}\right]}\nrightarrow
{\left[\begin{smallmatrix} \mu  &0\\
0&\nu
 \end{smallmatrix}\right]}
$  \emph{and}
${\left[\begin{smallmatrix}
\lambda  &0\\
0&\pm\lambda
 \end{smallmatrix}\right]}\nrightarrow
{\left[\begin{smallmatrix} 0 &1\\
\sigma &0
 \end{smallmatrix}\right]}
$ ($|\lambda|=|\mu |=|\nu |= 1$, $\mu
\ne\pm\nu$, $|\sigma |<1$).

Let there exist an arbitrarily small
perturbation
$A:=\left[\begin{smallmatrix}
\lambda  &0\\
0&\pm\lambda
 \end{smallmatrix}\right]+E$ of
 $\left[\begin{smallmatrix}
\lambda  &0\\
0&\pm\lambda
 \end{smallmatrix}\right]$ that is
 *congruent to
$B:={\left[\begin{smallmatrix} \mu  &0\\
0&\nu
 \end{smallmatrix}\right]}$ or
 $C:={\left[\begin{smallmatrix} 0 &1\\
\sigma &0
 \end{smallmatrix}\right]}$.
 Then $A^{-*}A:=(A^{-1})^*A$ is
similar to $B^{-*}B$ or $C^{-*}C$,
which is impossible since the
eigenvalues of $A^{-*}A$ are
arbitrarily close to $\bar \lambda
^{-1}\lambda =\lambda ^2$, whereas
$B^{-*}B=\diag(\mu^2,\nu ^2)$ and
$C^{-*}C= \diag(\sigma,\bar\sigma
^{-1})$.

\medskip

$\bullet$  ${\left[\begin{smallmatrix}
\lambda  &0\\
0&\lambda
 \end{smallmatrix}\right]}\nrightarrow
{\left[\begin{smallmatrix} 0 &\tau \\
\tau &\tau i
 \end{smallmatrix}\right]}
$ ($|\lambda  |=|\tau |=1$).

Let ${\left[\begin{smallmatrix}
\lambda  &0\\
0&\lambda
 \end{smallmatrix}\right]}
\to{\tau \left[\begin{smallmatrix} 0 &1\\
1&i
 \end{smallmatrix}\right]}
$; i.e., there exists an arbitrarily
small perturbation
$A:=\left[\begin{smallmatrix}
1 &0\\
0&1
 \end{smallmatrix}\right]+E$ of $
 \left[\begin{smallmatrix}
1 &0\\
0&1
 \end{smallmatrix}\right]$ that is
 *congruent to
$B:=\lambda^{-1}\tau
{\left[\begin{smallmatrix} 0 &1\\
1&i
 \end{smallmatrix}\right]}$.
This means that there exists a
nonsingular $S$ such that
\[S^*\left(\begin{bmatrix}
1 &0\\
0&1
 \end{bmatrix}+E\right)S=
 \lambda^{-1} \tau
 {\begin{bmatrix} 0 &1\\
1&i
 \end{bmatrix}}.\]
Equating the determinants of both
sides, we find that $(\lambda^{-1}
\tau)^2 =-1$, and so  $\lambda^{-1}
\tau =\pm i$. Then $\rank (B+B^*)=1$,
which is impossible since $A+A^*$ is
*congruent to $B+B^*$ and $\rank
(A+A^*)=2$.

\end{document}